\documentclass[preprints,article,accept,moreauthors,pdftex]{my-mdpi} 

\firstpage{1} 
\pubvolume{11}
\issuenum{4}
\articlenumber{178}
\doinum{10.3390/axioms11040178}
\pubyear{2022}
\copyrightyear{2022}
\externaleditor{Academic Editor: Chris Goodrich} 
\datereceived{21 March 2022} 
\daterevised{11 April 2022} 
\dateaccepted{13 April 2022} 
\datepublished{15 April 2022} 
\hreflink{https://doi.org/10.3390/\linebreak axioms11040178} 

\pdfoutput=1

\Title{Weighted Generalized Fractional Integration by Parts and the Euler--Lagrange Equation}

\TitleCitation{Weighted Generalized Fractional Integration by Parts and the Euler--Lagrange Equation}

\Author{Houssine Zine $^{1,\dagger}$\orcidA{},
El Mehdi Lotfi $^{2,\dagger}$\orcidB{}, 
Delfim F. M. Torres $^{1,}$*$^{,\dagger,}$\orcidC{} 
and Noura Yousfi $^{2,\dagger}$\orcidD{}}

\AuthorNames{Houssine Zine, El Mehdi Lotfi, Delfim F. M. Torres and Noura Yousfi}

\AuthorCitation{Zine, H.; Lotfi, E.M.; Torres, D.F.M.; Yousfi, N.}

\address{$^{1}$ \quad Center for Research and Development in Mathematics and Applications (CIDMA), 
Department of Mathematics, University of Aveiro, 3810-193 Aveiro, Portugal; 
zinehoussine@ua.pt\\
$^{2}$ \quad Laboratory of Analysis, Modeling and Simulation (LAMS),
Faculty of Sciences Ben M'sik, Hassan II University of Casablanca, 
P.O. Box 7955, Sidi Othman, Casablanca 20000, Morocco; 
\mbox{lotfiimehdi@gmail.com (E.M.L.);} nourayousfi.fsb@gmail.com (N.Y.)}

\corres{Correspondence: delfim@ua.pt}

\firstnote{These authors contributed equally to this work.}

\abstract{Integration by parts plays a crucial role in mathematical analysis, 
e.g., during the proof of necessary optimality conditions in the
calculus of variations and optimal control. Motivated by this fact, 
we construct a new, right-weighted generalized fractional derivative 
in the Riemann--Liouville sense with its associated integral for the 
recently introduced weighted generalized fractional derivative
with Mittag--Leffler kernel. We rewrite these operators equivalently 
in effective series, proving some interesting properties relating to
the left and the right fractional operators. These results 
permit us to obtain the corresponding integration by parts formula. 
With the new general formula, we obtain an appropriate weighted 
Euler--Lagrange equation for dynamic optimization,
extending those existing in the literature. We end with 
the application of an optimization variational problem
to the quantum mechanics framework.}



\begin{document}
	

\section{Introduction}
\label{sec:01}

\textls[-5]{In the last decade, fractional calculus played an important role
in the theoretical study of dynamical systems by showing significant 
results in many natural fields and engineering domains~\cite{MR4200078,lot1}. 
For this reason, mathematicians are paying more attention to the
generalization of several important formulas in the integral theory 
of Mathematical Analysis, namely, the Newton--Leibniz formula, 
the Green formula, and the Gauss and Stokes formulas~\cite{Stocks2,Stocks1}. 
Some are central tools that enable mathematicians to extend 
other theories, such as the integration by parts formula, 
Taylor's formula, the Euler--Lagrange equation, Gr\"{o}nwall's inequality, 
Lyapunov theorems and LaSalle's invariance principle~\cite{lasalle,Lyapunov}.}

Often, memory effects are fractionally modeled
with Riemann--Liouville and Caputo derivatives~\cite{MR4364786,Malinowska}. 
However, the fact that the Mittag--Leffler function is a generalization 
of the exponential function naturally gives rise to new definitions 
for fractional operators~\cite{MR4395549,MR4404431}. In 2020, Hattaf~\cite{hattaf} 
has proposed a new left-weighted generalized fractional derivative 
for both Caputo and Riemann--Liouville senses and their associated integral operator, 
see also~\cite{hattaf2}. Motivated by their applications in mechanics, 
where the introduction of the correct operator is needed~\cite{Malinowska,MR1438729}, 
here, we introduce the right-weighted generalized fractional derivative 
and its associated integral operator, proving their main properties
and, in particular, their integration by parts formula. 

It is worth emphasizing that integration by parts is of great interest 
in integral calculus and mathematical analysis. For example,
it represents a strong tool to develop the calculus of variations 
through the so-called Euler--Lagrange equation, which is the 
central result of dynamic optimization~\cite{Malinowska}. 
In recent years, the development of some theoretical practices using 
fractional derivatives has drawn the attention of several researchers. In 2012,
Almeida, Malinowska and Torres~\cite{tor} reviewed some recent results 
of fractional variational calculus and discussed the necessary optimality 
conditions of Euler--Lagrange type for functionals with a Lagrangian 
containing left and right Caputo derivatives. In 2017, Abdeljawad and Baleanu 
obtained an adequate integration by parts formula and the corresponding 
Euler--Lagrange equations using the nonlocal fractional derivative 
with Mittag--Leffler kernel. In 2019, Abdeljawad et al.~\cite{integ2}
developed a fractional integration by parts formula for 
Riemann--Liouville, Liouville--Caputo, Caputo--Fabrizio and Atangana--Baleanu 
fractional derivatives. In 2020, Zine and Torres~\cite{zine}
introduced a stochastic fractional calculus, and obtained 
a stochastic fractional Euler--Lagrange equation. 
Motivated by these works, particularly~\cite{integ2,integ1,tor,zine}, 
and with the help of our weighted generalized fundamental 
integration by parts formula, we extend the 
available Euler--Lagrange equations. 

The main purpose of our work is to compute a new integration by parts
formula for the weighted generalized fractional derivative and to discuss 
the associated necessary optimality conditions of Euler--Lagrange type.
To do this, we organize the paper as follows. In Section~\ref{sec:02}, 
we recall some necessary results from the literature. 
We proceed with \mbox{Section~\ref{sec:03},} introducing the right-weighted generalized 
fractional derivative and its associated integral and studying their well-posedness. 
Integration by parts is investigated in Section~\ref{sec:04}, followed
by Section~\ref{sec:05}, where the weighted generalized fractional Euler--Lagrange
equation is rigorously proved. We end with Section~\ref{sec:06}, illustrating
the obtained theoretical results with their application in the quantum mechanics framework.


\section{Preliminaries}
\label{sec:02}

\textls[-25]{In this section, we present some definitions and properties
from the fractional calculus literature, which will help us
to prove our main results. In the text, $f\in H^{1}(a,b)$ 
is a sufficiently smooth function on $[a,b]$
with $a, b \in \mathbb{R}$. In addition, 
we adopt the following~notations:}
\begin{gather*}
\phi(\alpha):=\dfrac{1-\alpha}{B(\alpha)},
\quad \psi(\alpha):=\dfrac{\alpha}{B(\alpha)},
\end{gather*}
where $0\leq\alpha< 1$ and $ B(\alpha)$ is a normalization function 
obeying $B(0)=B(1)=1$. In the paper, we denote 
$$
\mu_\alpha:=\dfrac{\alpha}{1-\alpha}.
$$

\begin{Lemma}[See~\cite{Samko}]
\label{Samko}
Let $ \alpha > 0$, $p\geq1$, $q\geq 1$ 
and $\dfrac{1}{p}+\dfrac{1}{q}\leq 1+\alpha$ 
($p\neq 1$ and $q\neq 1$ in the case  
$\dfrac{1}{p}+\dfrac{1}{q}= 1+\alpha$). 
If $f\in L_p(a,b)$ and $g\in L_q(a,b)$, then
\begin{equation*}
\int_{a}^{b}f(x)\,^{RL}_{a,1}I^{\alpha}g(x) dx
=\int_{a}^{b}g(x)\,^{RL}I^{\alpha}_{b,1}f(x) dx,
\end{equation*}
where {$^{RL}_{a,1}I^{\alpha}$} is the left standard Riemann--Liouville
fractional integral of order $\alpha$ given by
 \begin{equation}
^{RL}_{a,1}I^{\alpha}f(x)
=\dfrac{1}{\Gamma(\alpha)}
\int^{x}_{a}(x-s)^{\alpha-1}f(s)ds,
\quad x > a,
\end{equation}
and $^{RL}I^{\alpha}_{b,1}$ is the right standard Riemann--Liouville
fractional integral of order $\alpha$ given by
 \begin{equation}
^{RL}I^{\alpha}_{b,1}f(x)
=\dfrac{1}{\Gamma(\alpha)}
\int^{b}_{x}(s-x)^{\alpha-1}f(s)ds,
\quad x < b.
\end{equation}
\end{Lemma}

\begin{Definition}[See~\cite{hattaf}]
Let $0\leq\alpha< 1$ and $\beta>0$. 
The left-weighted generalized fractional derivative
of order $\alpha$ of function $f$, 
in the Riemann--Liouville sense, is defined by
\begin{equation}
\label{NGDR:g}
\,_{a,w}^{R}D^{\alpha,\beta}f(x)
= \dfrac{1}{\phi(\alpha)}\dfrac{1}{w(x)}\dfrac{d}{dx}\int^{x}_{a}(wf)(s)
E_\beta\left[ -\mu_{\alpha}(x-s)^{\beta}\right]ds,
\end{equation}
where $E_{\beta}$ denotes the Mittag--Leffler 
function of parameter $\beta$ defined by
\begin{equation}
\label{mitag}
E_{\beta}(z)=\sum_{j=0}^{\infty}\dfrac{z^{j}}{\Gamma(\beta j+1)},
\quad z \in \mathbb{C},
\end{equation}
and $w\in C^{1}([a,b])$ with $ w,w'>0$.
The corresponding fractional integral is defined by
\begin{equation}
\label{NGI:g}
_{a,w}I^{\alpha,\beta}f(x)
= \phi(\alpha)f(x)+\psi(\alpha)\,^{RL}_{a,w}I^{\beta}f(x),
\end{equation}
where $^{RL}_{a,w}I^{\beta}$ is the standard weighted Riemann--Liouville
fractional integral of order $\beta$ given by
 \begin{equation}
^{RL}_{a,w}I^{\beta}f(x)=\dfrac{1}{\Gamma(\beta)}
\dfrac{1}{w(x)}\int^{x}_{a}(x-s)^{\beta-1}w(s)f(s)ds,
\quad x > a.
\end{equation}
\end{Definition}


\section{Well-Posedness of the Right-Weighted Fractional Operators}
\label{sec:03}

We denote the right-weighted generalized fractional derivative of order $\alpha$ in the
Riemann--Liouville sense by $^{R}D^{\alpha,\beta}_{b,w}$, and we define this so that 
the following identity occurs: 
$$
Q\left(_{a,w}^{R}D^{\alpha,\beta}f\right)(x)
=\left(^{R}D^{\alpha,\beta}_{b,w}Qf\right)(x)
$$
with $Q$ being the \emph{reflection operator}, that is, 
$(Qf)(x) = f(a + b - x)$ with function $f$ defined on the interval $[a, b]$.

\begin{Definition}[right-weighted generalized fractional derivative]
\label{NGDR:d}
Let $0\leq\alpha< 1$ and $\beta>0$. The right-weighted generalized fractional derivative
of order $\alpha$ of function $f$, in the Riemann--Liouville sense, is defined by
\begin{equation}
^{R}D^{\alpha,\beta}_{b,w}f(x)
= \dfrac{-1}{\phi(\alpha)}\dfrac{1}{w(x)}\dfrac{d}{dx}\int^{b}_{x}(wf)(s)
E_\beta\left[ -\mu_{\alpha}(s-x)^{\beta}\right]ds,
\end{equation}
where $w\in C^{1}([a,b])$ with $w,w'>0$.
\end{Definition}

To properly define the new right-weighted fractional integral,
we need to solve the equation $^{R}D^{\alpha,\beta}_{b,w}f(x)=u(x)$. 
We have
$$
^{R}D^{\alpha,\beta}_{b,w}f(x)=\,^{R}D^{\alpha,\beta}_{b,w}QQf(x)
=Q\,_{a,w}^{R}D^{\alpha,\beta}Qf(x)=u(t).
$$

Then,
$$
_{a,w}^{R}D^{\alpha,\beta}Qf(x)=Qu(x)
$$
and thus,
$$
Qf(x)=\phi(\alpha)Qu(x)+\psi(\alpha)\,_{a,w}^{RL}I^{\beta}Qu(x)
=\phi(\alpha)Qu(x)+\psi(\alpha)\,Q^{RL}I^{\beta}_{b,w}u(x),
$$
where $^{RL}I^{\beta}_{b,w}$ is the right-weighted standard Riemann--Liouville
fractional integral of order $\beta$ given by
\begin{equation}
\label{eq:bs}
^{RL}I^{\beta}_{b,w}f(x)=\dfrac{1}{\Gamma(\beta)}
\dfrac{1}{w(x)}\int^{b}_{x}(s-x)^{\beta-1}w(s)f(s)ds,
\quad x < b.
\end{equation}

Applying $Q$ to both sides of \eqref{eq:bs}, we obtain 
$$
f(t)=\phi(\alpha)u(x)+\psi(\alpha)\,^{RL}I^{\beta}_{b,w}u(x).
$$

Moreover,
\begin{eqnarray*}
_{a,w}I^{\alpha,\beta}Qf(x) 
&=& \phi(\alpha)Qf(x)+\psi(\alpha)\,^{RL}_{a,w}I^{\beta}Qf(x)\\
&=& \phi(\alpha)Qf(x)+\psi(\alpha)\,Q^{RL}I^{\beta}_{b,w}f(x) \\
&=& Q\big[\phi(\alpha)f(x)+\psi(\alpha)\,^{RL}I^{\beta}_{b,w}f(x)\big].
\end{eqnarray*}

We are now in the position to introduce the concept 
of the right-weighted generalized fractional integral.

\begin{Definition}[right-weighted generalized fractional integral]
Let $0\leq\alpha< 1$ and $\beta>0$.  The right-weighted generalized fractional integral
of order $\alpha$ of function $f$ is given by
\begin{equation}
\label{NGI:d}
I^{\alpha,\beta}_{b,w}f(x)
= \phi(\alpha)f(x)+\psi(\alpha)\,^{RL}I^{\beta}_{b,w}f(x),
\end{equation}
where $ w\in C^{1}([a,b])$ with $ w,w'>0$.
\end{Definition}

Our next result provides a series representation
to the left- and right-weighted generalized fractional derivatives.

\begin{Theorem}
\label{serie}
Let $0\leq\alpha< 1$ and $\beta>0$. The left- and 
right-weighted generalized fractional derivatives
of order $\alpha$ of function $f$ 
can be written, respectively, as
\begin{equation}
\label{seriea}
\,_{a,w}^{R}D^{\alpha,\beta}f(x)
=\dfrac{1}{\phi(\alpha)}\sum_{j=0}^{\infty}(-\mu_\alpha)^j
\,^{RL}_{a,w}I^{\beta j}f(x)
\end{equation}
and
\begin{equation}
\label{serieb}
\,^{R}D^{\alpha,\beta}_{b,w}f(x)
=\dfrac{-1}{\phi(\alpha)}\sum_{j=0}^{\infty}(-\mu_\alpha)^j
\,^{RL}_{b,w}I^{\beta j}f(x).
\end{equation}
\end{Theorem}

\begin{proof}
The Mittag--Leffler function $E_\beta(x)$ is an entire series of $x$. 
Since the series (\ref{mitag}) locally and uniformly  converges
in the whole complex plane, the left-weighted generalized 
fractional derivative can be rewritten as
\begin{eqnarray*}
\,_{a,w}^{R}D^{\alpha,\beta}f(x)
&=&  \dfrac{1}{\phi(\alpha)}\dfrac{1}{w(x)}\dfrac{d}{dx}\int^{x}_{a}(wf)(s)
\sum_{j=0}^{\infty}(-\mu_\alpha)^j\dfrac{(x-s)^{\beta j}}{\Gamma(\beta j+1)}ds \\
&=& \dfrac{1}{\phi(\alpha)}\dfrac{1}{w(x)} \sum_{j=0}^{\infty} 
(-\mu_\alpha)^j \dfrac{1}{\Gamma(\beta j+1)}
\dfrac{d}{dx} \int^{x}_{a}(wf)(s)(x-s)^{\beta j} ds\\
&=&\dfrac{1}{\phi(\alpha)}\dfrac{1}{w(x)} \sum_{j=0}^{\infty} 
(-\mu_\alpha)^j \dfrac{1}{\Gamma(\beta j)}
\int^{x}_{a}(wf)(s)(x-s)^{\beta j-1} ds\\
&=& \dfrac{1}{\phi(\alpha)}\sum_{j=0}^{\infty}
(-\mu_\alpha)^j\big(\,^{RL}_{a,w}I^{\beta j}f(x)\big).
\end{eqnarray*}

From Definition~\ref{NGDR:d}, and using the same steps that were used before, 
one can easily rewrite the new right-weighted generalized 
fractional derivative as equality (\ref{serieb}).
The proof of \eqref{seriea} is similar.
\end{proof}

\begin{Theorem}
\label{NL}
Let $0\leq\alpha< 1$ and $\beta>0$. The left- and right-weighted 
generalized fractional derivative and their associated integrals
satisfy the following formulas:
\begin{equation}
\label{NLa}
{_{a,w}I^{\alpha,\beta}}\big({_{a,w}^{R}D^{\alpha,\beta}}f\big)(x)
=_{a,w}^{R}D^{\alpha,\beta}\big(\,_{a,w}I^{\alpha,\beta}f\big)(x)=f(x)
\end{equation}
and
\begin{equation}
\label{NLb}
I^{\alpha,\beta}_{b,w}\big(^{R}D_{b,w}^{\alpha,\beta}f\big)(x)
=^{R}D_{b,w}^{\alpha,\beta}\big(\,^{\alpha,\beta}I_{b,w}f\big)(x)=-f(x).
\end{equation}
\end{Theorem}

\begin{proof}
We note that
\begin{equation*}
\begin{split}
_{a,w}I^{\alpha,\beta}\big(_{a,w}^{R}D^{\alpha,\beta}f\big)(x) 
&= \phi(\alpha)\big(_{a,w}^{R}D^{\alpha,\beta}f\big)(x)
+\psi(\alpha)^{RL}_{a,w}I^{\beta}\big(_{a,w}^{R}D^{\alpha,\beta}f\big)(x)\\
&=  \sum_{j=0}^{\infty}(-\mu_\alpha)^j\,^{RL}_{a,w}I^{\beta j}f(x)
+\mu_\alpha\,^{RL}_{a,w}I^{\beta}\big(\sum_{j=0}^{\infty}
(-\mu_\alpha)^j\,^{RL}_{a,w}I^{\beta j}f\big)(x) \\
&= \sum_{j=0}^{\infty}(-\mu_\alpha)^j\,^{RL}_{a,w}I^{\beta j}f(x)
-\sum_{j=0}^{\infty}(-\mu_\alpha)^{j+1}\,^{RL}_{a,w}I^{\beta + \beta j}f(x)\\
&=  f(t).
\end{split}
\end{equation*}

Then,
\begin{equation*}
\begin{split}
_{a,w}^{R}D^{\alpha,\beta}\big(\,_{a,w}I^{\alpha,\beta}f\big)(x) 
&= \dfrac{1}{\phi(\alpha)}\sum_{j=0}^{\infty}
(-\mu_\alpha)^j\,^{RL}_{a,w}I^{\beta j}\big(\,_{a,w}I^{\alpha,\beta}f\big)(x),\\
&= \dfrac{1}{\phi(\alpha)}\sum_{j=0}^{\infty}
(-\mu_\alpha)^j\,^{RL}_{a,w}I^{\beta j}\bigg[\phi(\alpha)f(x)
+\psi(\alpha)\,^{RL}_{a,w}I^{\beta}f(x)\bigg] \\
&=  \sum_{j=0}^{\infty}(-\mu_\alpha)^j\,^{RL}_{a,w}I^{\beta j}f(x)
+\mu_\alpha\sum_{j=0}^{\infty}(-\mu_\alpha)^j\,^{RL}_{a,w}I^{\beta j+\beta}f(x)\\
&=  \sum_{j=0}^{\infty}(-\mu_\alpha)^j\,^{RL}_{a,w}I^{\beta j}f(x)-\sum_{j=0}^{\infty}
(-\mu_\alpha)^{j+1}\,^{RL}_{a,w}I^{\beta j+\beta}f(x) \\
&= f(x)
\end{split}
\end{equation*}
and equality (\ref{NLa}) holds true. The proof 
of equality (\ref{NLb}) is similar.
\end{proof}


\section{Integration by Parts}
\label{sec:04}

Our formulas of integration by parts 
are proved in suitable function spaces.

\begin{Definition}[See~\cite{kilbas}]
For $\alpha>0$, $\beta>0$ and $1\leq p\leq \infty$, 
the following function spaces are~defined:
$$
_{a,w}I^{\alpha,\beta}(L_p)
:=\left\{f:f=\,_{a,w}I^{\alpha,\beta}(\eta), 
\eta \in L_p(a,b)\right\}
$$
and
$$
I_{b,w}^{\alpha,\beta}(L_p)
:=\left\{f:f=I^{\alpha,\beta}_{b,w}(\theta), 
\theta \in L_p(a,b)\right\}.
$$
\end{Definition}

\begin{Theorem}[integration by parts without the weighted function]
\label{IBPWWF}
Let $0 \leq \alpha < 1$, $\beta >0$, $p\geq1$, $q\geq 1$ 
and $\dfrac{1}{p}+\dfrac{1}{q}\leq 1+\alpha$ 
{($p\neq 1$ and $q\neq 1$ in the case 
$\dfrac{1}{p}+\dfrac{1}{q} = 1+\alpha$)}.
\begin{itemize}
\item If $ f\in L_p(a,b)$  and $g\in L_q(a,b)$, then
\begin{equation}
\label{I1}
\int_{a}^{b}f(x)(_{a,1}I^{\alpha,\beta}g)(x)dx
= \int_{a}^{b} g(x)(I^{\alpha,\beta}_{b,1}f)(x)dx.
\end{equation}

\item If $ f\in I^{\alpha,\beta}_{b,w}(L_p)$  
and $g\in \,_{a,w}I^{\alpha,\beta}(L_q)$, then
\begin{equation}
\label{D1}
\int_{a}^{b}f(x)(_{a,1}^RD^{\alpha,\beta}g)(x)dx
= \int_{a}^{b} g(x)(^RD^{\alpha,\beta}_{b,1}f)(x)dx.
\end{equation}
\end{itemize}
\end{Theorem}

\begin{proof}
First, we prove equality (\ref{I1}). Since,
\begin{equation*}
\begin{split}
\int_{a}^{b}f(x)(_{a,1}I^{\alpha,\beta}g)(x) dx 
&= \int_{a}^{b}f(x)\left[\phi(\alpha)g(x)
+\psi(\alpha)\,^{RL}_{a,1}I^{\beta}g(x)\right]\\
&=  \phi(\alpha)\int_{a}^{b}f(x)g(x)dx
+\psi(\alpha)\int_{a}^{b}f(x)\,^{RL}_{a,1}I^{\beta}g(x) dx,
\end{split}
\end{equation*}
it follows from Lemma~\ref{Samko} that
\begin{equation*}
\begin{split}
\int_{a}^{b}f(x)(_{a,1}I^{\alpha,\beta}g)(x) dx 
&=  \phi(\alpha)\int_{a}^{b}f(x)g(x)dx
+\psi(\alpha)\int_{a}^{b}g(x)\,^{RL}I^{\beta}_{b,1}f(x) dx\\
&= \int_{a}^{b}g(x)\left[\phi(\alpha)f(x)
+\psi(\alpha)\,^{RL}I^{\beta}_{b,1}f(x)\right]\\
&= \int_{a}^{b}g(x)(I^{\alpha,\beta}_{b,1}f)(x) dx.
\end{split}
\end{equation*}

Now, we prove (\ref{D1}):
\begin{equation*}
\begin{split}
\int_{a}^{b}f(x)(_{a,1}^RD^{\alpha,\beta}g)(x)dx 
&= \int_{a}^{b}I^{\alpha,\beta}_{b,1}\theta(x)
\left(_{a,1}^RD^{\alpha,\beta}(\,_{a,1}I^{\alpha,\beta}\eta)\right)(x)dx \\
&= \int_{a}^{b}\eta(x)I^{\alpha,\beta}_{b,1}\theta(x)dx  
\ \ \big(\text{from Theorem~\ref{NL}}\big)\\
&= \int_{a}^{b}\theta(x)\,_{a,1}I^{\alpha,\beta}\eta(x)dx 
\ \ \big(\text{from equality} \ (\ref{I1})\big)\\
&=\int_{a}^{b}g(x)(\,^RD^{\alpha,\beta}_{b,1}f)(x)dx 
\ \ \big(\text{from Theorem~\ref{NL}}\big).
\end{split}
\end{equation*}

The proof is complete.
\end{proof}

\begin{Theorem}[weighted generalized integration by parts]
\label{GIBP}
Let $0 \leq \alpha < 1$, $\beta >0$, $p\geq1$, $q\geq 1$ and 
$\dfrac{1}{p}+\dfrac{1}{q}\leq 1+\alpha $ ($p\neq 1$ and $q\neq 1$ 
in the case $\dfrac{1}{p}+\dfrac{1}{q} = 1+\alpha$). 
If $f\in L_p(a,b)$ and $g\in L_q(a,b)$, then
\begin{eqnarray}
\int_{a}^{b}f(x)(_{a,w}I^{\alpha,\beta}g)(x)dx
&=& \int_{a}^{b} w(x)^2g(x)\left(I^{\alpha,\beta}_{b,w}
\left(\dfrac{f}{w^2}\right)\right)(x)dx,\label{Iaw}\\
\int_{a}^{b}f(x)(\,^R_{a,w}D^{\alpha,\beta}g)(x)dx
&=& \int_{a}^{b} w(x)^2g(x)\left(\,^RD^{\alpha,\beta}_{b,w}
\left(\dfrac{f}{w^2}\right)\right)(x)dx.\label{RDaw}
\end{eqnarray}
\end{Theorem}

\begin{proof}
We have
\begin{equation*}
\begin{split}
\int_{a}^{b}f(x)(_{a,w}I^{\alpha,\beta}g)(x)dx 
&= \int_{a}^{b}w(x)\dfrac{f(x)}{w(x)}
\left(\,_{a,w}I^{\alpha,\beta}\left(
\dfrac{gw}{w}\right)\right)(x)dx \\
&=  \int_{a}^{b}\dfrac{f(x)}{w(x)}
\left(\,_{a,1}I^{\alpha,\beta}\left(gw\right)\right)(x)dx\\
&=  \int_{a}^{b}w(x)g(x)\left(I^{\alpha,\beta}_{b,1}\left(
\dfrac{f}{w}\right)\right)(x)dx 
\ \text{\big (from Theorem~\ref{IBPWWF}\big)}\\
&=  \int_{a}^{b}g(x)w(x)^2\left(I^{\alpha,\beta}_{b,w}\left(
\dfrac{f}{w^2}\right)\right)(x)dx.
\end{split}
\end{equation*}

Therefore, equality (\ref{Iaw}) is true. Similarly, we have
\begin{equation*}
\begin{split}
\int_{a}^{b}f(x)(_{a,w}^RD^{\alpha,\beta}g)(x)dx 
&= \int_{a}^{b}w(x)\dfrac{f(x)}{w(x)}
\left(\,_{a,w}^RD^{\alpha,\beta}\left(\dfrac{gw}{w}\right)\right)(x)dx \\
&=  \int_{a}^{b}\dfrac{f(x)}{w(x)}
\left(\,_{a,1}^RD^{\alpha,\beta}\left(gw\right)\right)(x)dx\\
&=  \int_{a}^{b}w(x)g(x)\left(D^{\alpha,\beta}_{b,1}
\left(\dfrac{f}{w}\right)\right)(x)dx 
\ \text{\big (from Theorem~\ref{IBPWWF}\big)}\\
&=  \int_{a}^{b}g(x)w(x)^2\left(D^{\alpha,\beta}_{b,w}
\left(\dfrac{f}{w^2}\right)\right)(x)dx
\end{split}
\end{equation*}
and equality (\ref{RDaw}) holds.
\end{proof}

\begin{Remark}
When $w(t)=1$ and $\alpha=\beta$, then we can obtain 
from our Theorem~\ref{GIBP} the integration by parts 
formula \cite{integ1} associated with Atangana--Baleanu 
derivatives:
\begin{eqnarray*}
\int_{a}^{b}f(x)\left(_{a}^{AB}I^{\alpha}g\right)(x)dx
&=& \int_{a}^{b} g(x)\left(^{AB}I^{\alpha}_{b} f\right)(x)dx,\\
\int_{a}^{b}f(x)\left(^{ABR}_{a}D^{\alpha}g\right)(x)dx
&=& \int_{a}^{b} g(x)\left(^{ABR}D^{\alpha}_{b}f\right)(x)dx.
\end{eqnarray*}
\end{Remark}

From (\ref{Iaw}) and (\ref{RDaw}),
we obtain the following consequence.

\begin{Corollary}
Let $0 \leq \alpha < 1$, $\beta >0$, $p\geq1$, $q\geq 1$ and 
$\dfrac{1}{p}+\dfrac{1}{q}\leq 1+\alpha $ ($p\neq 1$ and $q\neq 1$ 
in the case $\dfrac{1}{p}+\dfrac{1}{q} = 1+\alpha$). 
If $f\in L_p(a,b)$ and $g\in L_q(a,b)$, then
\begin{equation*}
\begin{split}
\int_{a}^{b}f(x)\left(I^{\alpha,\beta}_{b,w}g\right)(x)dx
&= \int_{a}^{b} w(x)^2g(x)\left(\,_{a,w}I^{\alpha,\beta}
\left(\dfrac{f}{w^2}\right)\right)(x)dx,\\
\int_{a}^{b}f(x)\left(\,^RD^{\alpha,\beta}_{b,w}g\right)(x)dx
&= \int_{a}^{b} w(x)^2g(x)\left(\,_{a,w}^RD^{\alpha,\beta}
\left(\dfrac{f}{w^2}\right)\right)(x)dx.
\end{split}
\end{equation*}
\end{Corollary}

For a symmetric view of weighted generalized integration by parts, 
we propose the following corollary of Theorem~\ref{GIBP}.

\begin{Corollary}
Let $0 \leq \alpha < 1$, $\beta >0$, $p\geq1$, $q\geq 1$ 
and $\dfrac{1}{p}+\dfrac{1}{q}\leq 1+\alpha $ ($p\neq 1$ 
and $q\neq 1$ in the case $\dfrac{1}{p}+\dfrac{1}{q} = 1+\alpha$). 
If $f\in L_p(a,b)$ and $g\in L_q(a,b)$, then
\begin{eqnarray}
\int_{a}^{b}w(x)f(x)\left(\,_{a,w}I^{\alpha,\beta}\dfrac{g}{w}\right)(x)dx
&=& \int_{a}^{b} w(x)g(x)\left(I^{\alpha,\beta}_{b,w}
\dfrac{f}{w}\right)(x)dx,\label{Iaws}\\
\int_{a}^{b}w(x)f(x)\left(\,_{a,w}^RD^{\alpha,\beta}\dfrac{g}{w}\right)(x)dx
&=& \int_{a}^{b} w(x)g(x)\left(\,^RD^{\alpha,\beta}_{b,w}
\dfrac{f}{w}\right)(x)dx.\label{RDaws}
\end{eqnarray}
\end{Corollary}


\section{The Weighted Generalized Fractional Euler--Lagrange Equation}
\label{sec:05}

Let us denote by $AC(I\rightarrow \mathbb{R})$ the set of 
absolutely continuous functions $X$, where $I=[a,b]$, such that 
the left and right Riemann--Liouville-weighted generalized 
fractional derivatives of $X$ exist, endowed with the norm
$$
\Vert X \Vert = \sup_{t\in I}\left(\mid X(t) \mid
+\mid  \,_{a,w}^{RL}D^{\alpha,\beta}X(t) \mid
+\mid  ^{RL}D^{\alpha,\beta}_{b,w}X(t) \mid\right).
$$

Let $L \in C^{1}(I\times  \mathbb{R} \times\mathbb{R}
\times \mathbb{R}\rightarrow\mathbb{R})$
and consider the following minimization problem:
\begin{equation}
\label{eq:F}
J[X]=\left(\int_a^b
L\left( t,X(t), ^{RL}D^{\alpha,\beta}_{a,w}X(t),
^{RL}D^{\alpha,\beta}_{b,w}X(t)\right)dt\right)
\longrightarrow \min
\end{equation}
subject to the boundary conditions
\begin{equation}
\label{eq:BC}
X(a)=X_a, \quad X(b)=X_b.
\end{equation}

Under appropriate general conditions, one can prove that the minimum of $J[\cdot]$ 
exists~\cite{MR3200762}. Here, we are interested in showing the usefulness 
of our Theorem~\ref{GIBP} to prove the necessary optimality conditions 
for problem \eqref{eq:F} and \eqref{eq:BC}.
With the help of weighted generalized fractional integration by parts,
we obtain the following Euler--Lagrange necessary optimality condition 
for the fundamental weighted generalized fractional problem 
of the calculus of variations \eqref{eq:F} and \eqref{eq:BC}. 

\begin{Theorem}[the weighted generalized fractional Euler--Lagrange equation]
\label{thm:SFE-Leq}
If $L \in C^{1}(I \times \mathbb{R} \times \mathbb{R}
\times \mathbb{R} \rightarrow \mathbb{R})$
and $X\in AC([a,b]\rightarrow \mathbb{R})$ is a minimizer
of \eqref{eq:F} subject to the fixed end points \eqref{eq:BC};
then, $X$ satisfies the following weighted generalized 
fractional Euler--Lagrange equation:
$$
\partial_{2}L + w(t)^2  \,^{R}D^{\alpha,\beta}_{b,w}
\left(\frac{\partial_{3}L}{w(t)^2}\right)
+w(t)^2  \,^{R}_{a,w}D^{\alpha,\beta}
\left(\dfrac{\partial_{4}L}{w(t)^2}\right)=0,
$$
where $\partial_i L$ denotes the partial derivative of the Lagrangian
$L$ with respect to its $i$th argument evaluated at
$\left(t,X(t),^{RL}D^{\alpha,\beta}_{a,w}X(t),^{RL}D^{\alpha,\beta}_{b,w}X(t)\right)$.
\end{Theorem}

\begin{proof}
Let
$J[X]=\displaystyle \int_a^b L\left(t,X(t),^{R}_{a,w}D^{\alpha,\beta}X(t),
^{R}D^{\alpha,\beta}_{b,w} X(t)\right)dt$ and
assume that $X^*$ is the optimal solution
of problem \eqref{eq:F} and \eqref{eq:BC}. Set
$$
X=X^*+\varepsilon \eta,
$$
where $\eta, X\in AC([a,b]\rightarrow \mathbb{R})$
and $\varepsilon$ is a small, real parameter. By linearity
of the weighted generalized fractional derivative, we obtain
$$
^{R}_{a,w}D^{\alpha,\beta}X(t)=\,^{R}_{a,w}D^{\alpha,\beta}X^*
+\varepsilon\left( \, ^{R}_{a,w}D^{\alpha,\beta} \eta(t) \right)
$$
and
$$
^{R}D^{\alpha,\beta}_{b,w}X(t)
=\,^{R}D^{\alpha,\beta}_{b,w}X^*
+\varepsilon\left( \,^{R}D^{\alpha,\beta}_{b,w} \eta(t)\right).
$$

Now, consider the following function:
\begin{multline*}
J(\varepsilon)=\int_a^b L\left(t,X^*(t)
+\varepsilon \eta(t), \ ^{R}_{a,w}D^{\alpha,\beta} X^*(t)
+\varepsilon\left( ^{R}_{a,w}D^{\alpha,\beta} \eta(t)\right),\right.\\
\left.^{R}D^{\alpha,\beta}_{b,w} X^*(t)
+\varepsilon\left( ^{R}D_{b,w}^{\alpha,\beta} \eta(t)\right)\right) dt.
\end{multline*}

Fermat's theorem asserts that 
$\left.\dfrac{d}{d\varepsilon}J(\varepsilon)\right|_{\varepsilon=0} = 0$
and we deduce, by the chain rule, that
$$
\int_a^b\left(\partial_{2}L \cdot \eta
+\partial_{3}L \cdot ^{R}_{a,w}D^{\alpha,\beta}\eta
+\partial_{4}L \cdot ^{R}D^{\alpha,\beta}_{b,w}\eta\right)dt =0.
$$

Using Theorem~\ref{GIBP}
of weighted fractional integration by parts, we obtain that
$$
\int_a^b\left(\partial_{2}L \cdot \eta+ w(t)^2 \cdot 
\eta \cdot \,^{R}D^{\alpha,\beta}_{b,w}
\left(\frac{\partial_{3}L}{w(t)^2}\right)
+w(t)^2 \cdot \eta \cdot \,^{R}_{a,w}D^{\alpha,\beta}
\left(\dfrac{\partial_{4}L}{w(t)^2}\right)\right) dt =0.
$$

The result follows by the fundamental theorem 
of the calculus of variations.
\end{proof}


\section{An Application}
\label{sec:06}

\textls[-25]{Let us consider the weighted generalized 
fractional variational problem
\eqref{eq:F} and~\eqref{eq:BC}~with}
\begin{multline*}
L\left(t,X(t),^{R}_{a,w}D^{\alpha,\beta}X(t),
^{R}D^{\alpha,\beta}_{b,w}X(t)\right)\\
=\frac{1}{2}\left(\frac{1}{2}m
\mid ^{R}_{a,w}D^{\alpha,\beta}X(t) \mid ^2
+\frac{1}{2}m\mid ^{R}D^{\alpha,\beta}_{b,w}X(t)\mid ^2\right)-V(X(t)),
\end{multline*}
where $X$ is an absolutely continuous function on $[a,b]$ and 
$V$ maps $C^1(I \rightarrow \mathbb{R})$ to $\mathbb{R}$.
Note~that
$$
\frac{1}{2}\left(\frac{1}{2}m\mid ^{R}_{a,w}D^{\alpha,\beta}X(t) 
\mid ^2 +\frac{1}{2}m \mid^{R}D^{\alpha,\beta}_{b,w}X(t)\mid ^2\right)
$$
can be viewed as a weighted generalized kinetic energy in
the quantum mechanics framework. By applying our
Theorem~\ref{thm:SFE-Leq} to the current variational
problem, we obtain that
\begin{equation}
\label{eq:Newton:dyn:law}
\frac{1}{2}m \left[ w(t)^2\,^{R}D^{\alpha,\beta}_{b,w}
\left(\frac{\,^{R}_{a,w}D^{\alpha,\beta}X(t)}{w(t)^2}\right)
+w(t)^2\,^{R}_{a,w}D^{\alpha,\beta}\left(
\frac{\,^{R}D_{b,w}^{\alpha,\beta}X(t)}{w(t)^2}\right)\right]
=V'(X(t)),
\end{equation}
where $V'$ is the derivative of the potential energy of the system.
We observe that relation~\eqref{eq:Newton:dyn:law} generalizes 
Newton's dynamical law $m \ddot{X}(t)=V'(X(t))$.


\section{Conclusions}

In this work, some definitions and properties of a recent class 
of fractional operators defined by general integral operators, 
with and without singular kernels, are recalled. A new definition 
of a right-weighted generalized fractional operator in the  
Riemann--Liouville sense is then posed, serving as a 
prerequisite for the establishment of a new weighted generalized 
integration by parts formula, which shows a duality relation 
with the existing left weighted generalized fractional operator 
in the  Riemann--Liouville--Hattaf sense~\cite{hattaf}. 
In the context of the fractional calculus of variations,
we have investigated weighted generalized  Euler--Lagrange equations, 
which were then used to produce an effective application in 
the quantum mechanics setting, after a proper definition of kinetic energy. 


\vspace{6pt}
\authorcontributions{Conceptualization, H.Z., E.M.L., D.F.M.T. and N.Y.; 
validation, H.Z., E.M.L., D.F.M.T. and N.Y.; 
formal analysis, H.Z., E.M.L., D.F.M.T. and N.Y.; 
investigation, H.Z., E.M.L., D.F.M.T. and N.Y.; 
writing---original draft preparation, H.Z., E.M.L., D.F.M.T. and N.Y.; 
writing---review and editing, H.Z., E.M.L., D.F.M.T. and N.Y.; supervision, D.F.M.T. 
All authors have read and agreed to the published version of the manuscript.}

\funding{This research was funded by Funda\c{c}\~{a}o para a Ci\^{e}ncia e a Tecnologia (FCT)
grant number UIDB/04106/2020 (CIDMA).}

\dataavailability{Not applicable.} 

\acknowledgments{The authors would like to express their gratitude 
to two anonymous reviewers, for their constructive comments 
and suggestions, which helped them to enrich the paper.}

\conflictsofinterest{The authors declare no conflict of interest.
The funders had no role in the design of the study; in the collection, 
analyses, or interpretation of data; in the writing of the manuscript, 
or in the decision to publish the~results.} 

\end{paracol}


\reftitle{References}



\end{document}